\documentclass{amsart}
\usepackage{amscd} 
\usepackage{amsfonts} 
\usepackage{amssymb} 
\usepackage{latexsym} 

\newcommand{\ncm}{\newcommand}

\newtheorem{theorem}{Theorem}[section]
\newtheorem{prop}[theorem]{Proposition}
\newtheorem{lemma}[theorem]{Lemma}
\newtheorem{cor}[theorem]{Corollary}
\newtheorem{lem&def}[theorem]{Lemma \& Definition}
\newtheorem{definition}[theorem]{Definition}

\ncm{\End}{\mbox{\rm End}\,}

\def\Hom{\mbox{\rm Hom}\,}

\def\id{\mbox{\rm id}}

\def\to{\rightarrow}

\def\o{\otimes}    
\def\bra{\langle}
\def\ket{\rangle}

\ncm{\rarr}[1]{\stackrel{#1}{\longrightarrow}}
\ncm{\larr}[1]{\stackrel{#1}{\longleftarrow}}
\def\cop{\Delta}

\def\eps{\varepsilon}

\def\du1{\hat 1}

\def\-1{_{(-1)}}
\def\0{_{(0)}}
\def\1{_{(1)}}
\def\2{_{(2)}}
\def\3{_{(3)}}

\def\du1{\hat 1}

\begin{document}

\title[An approach to quasi-Hopf algebras]{An approach to quasi-Hopf algebras via Frobenius
coordinates}
\author{Lars Kadison} 
\address{Matematiska Institutionen \\ G{\" o}teborg 
University \\ 
S-412 96 G{\" o}teborg, Sweden} 
\email{lkadison@c2i.net} 
\thanks{}
\date{} 
\begin{abstract}
We study  quasi-Hopf algebras and their subobjects
 over certain commutative rings from the point
of view of Frobenius algebras. We introduce a 
type of Radford formula involving 
an anti-automorphism and the Nakayama automorphism of a Frobenius
algebra, then view several results in quantum algebras from
this vantage-point. In addition,  
 separability and strong separability
of quasi-Hopf algebras are studied as Frobenius algebras.  
\end{abstract} 
\maketitle

\section*{Dedicated to A.A.\ Stolin on his fiftieth birthday}

In \cite{D}, Drinfel'd introduces quasi-Hopf algebras over a commutative ground ring,
and works out the fundamentals of this theory of quasi-bialgebra with antipode.  
From a categorical point of view,
modules over a quasi-bialgebra form a monoidal category with a nontrivial associativity
constraint.  Conjugating the comultiplication of a bialgebra by 
a gauge element produces
nontrivial examples of quasi-bialgebras. Quasi-bialgebras then differ from bialgebras by  being
 only coassociative up to conjugation by a three-cocycle; cf.\ Eqs.~(\ref{eq: quasi-coassoc})
and~(\ref{eq: 3-cocycle}). Several applications are made by Drinfel'd to the Knizhnik-Zamolodchikov
system of p.d.e.'s and to Reshetikhin's method for obtaining knot invariants. 

Hausser and Nill have shown in \cite{HN} that finite-dimensional quasi-Hopf algebras over fields are
Frobenius algebras.
We 
would like to return in this paper to the  general
 commutative ground ring  for 
 quasi-Hopf algebras as much as possible while retaining aspects of Frobenius algebras.   
In the preliminaries, we first show that quasi-Hopf algebras over a commutative ring $k$ with trivial
Picard group are Frobenius $k$-algebras by sketching the direct approach of 
Bulacu-Caenepeel \cite{BC} to the isomorphism $\int^{\ell}_H \o H^* \cong H$
 of a quasi-Hopf algebra $H$,
its dual $H^*$ and its space of left integrals $\int^{\ell}_H$
 via the Van Daele-Panaite-Van Oystaeyen projection $P : H \to \int^{\ell}_H$. 
Somewhat more generally, we introduce
QFH-algebras, which are quasi-Hopf algebras over commutative rings that are Frobenius
algebras. We then study a Frobenius coordinate system derived from 
 \cite{BC}, transform it to the Frobenius system introduced in Hausser-Nill \cite{HN}
 and make various deductions starting from a type of Radford formula for an anti-automorphism of a Frobenius algebra
(Lemma~\ref{lem-preRadford}). First, 
 a simplified proof and extension of the Hausser-Nill-Radford formula for the fourth
power of the antipode is provided for QFH-algebras (Theorem~\ref{th-HN}). Second, a quasi-Hopf subalgebra, 
stable under an antipode  of $H$, is a $\beta$-Frobenius extension (Theorem~\ref{th-beta}). 
In  section~2 we make a study when a quasi-Hopf algebra is separable
or strongly separable in the sense of Kanzaki.


\section{Preliminaries on quasi-Hopf algebras}
In this section we review the basics of the Bulacu-Caenepeel approach to quasi-Hopf algebras
with  small changes in notational conventions, generality and closer attention to Frobenius
systems. 

Let $k$ be a commutative ground ring.  All unlabelled tensors and Hom's are over
$k$.  Let $H$ be a finitely generated projective $k$-module and
$k$-algebra. All unlabeled identity maps and unity elements are on or in $H$.
Note that such an algebra is  Dedekind-finite in that $xy = 1$ if and only
if $yx = 1$.    We let $H^*$ denote $\Hom(H,k)$, which has the natural $H$-bimodule structure
defined as usual by $\bra h \rightharpoonup h^* \leftharpoonup k | x \ket = \bra h^* | kxh \ket$
(this notation is consistent with customary Hopf algebra notation). 

Following Drinfeld, we say that $H$ is a \textit{quasi-bialgebra} if it admits additional
structure $(H, \cop, \eps, \Phi)$ where algebra homomorphism $\cop: H \to H \otimes H$ 
is a possibly noncoassociative coproduct
with counit augmentation $\eps: H \to k$ satisfying the ordinary counit laws:
 $$(\eps \otimes \id)\cop = \id = (\id \otimes \eps)\cop.$$
$\Phi$ is an invertible element in $H \otimes H \otimes H$ denoted by 
$$\Phi = X^1 \otimes X^2 \otimes X^3 = Y^1 \otimes Y^2 \otimes Y^3 = \cdots $$
with inverse denoted by
$$\Phi^{-1} = x^1 \otimes x^2 \otimes x^3 = y^1 \otimes y^2 \otimes y^3 = \cdots $$
where we suppress a possible summation and change capital or lowercase letters for each occurence of $\Phi$ or $\Phi^{-1}$,
respectively, in the same side of an equation.
$\Phi$ controls the noncoassocativity of the coproduct on $H$ as follows:
\begin{equation}
\label{eq: quasi-coassoc}
(\id \otimes \cop)(\cop(a)) = \Phi (\cop \otimes \id)(\cop(a)) \Phi^{-1}
\end{equation}
which in generalized Sweedler notation is equivalent to 
$$a\1 X^1 \otimes a_{(2,1)} X^2 \otimes a_{(2,2)} X^3 = X^1 a_{(1,1)} \otimes X^2 a_{(1,2)} \otimes X^3 a\2. $$
Moreover $\Phi$ must satisfy normalized $3$-cocycle equations given by (an equation in $H \otimes H \otimes H \otimes H$
and another in $H \otimes H$): 
\begin{equation}
\label{eq: 3-cocycle}
(1 \otimes \Phi)(\id \otimes \cop \otimes \id)(\Phi)(\Phi \otimes 1) = 
(\id \otimes \id \otimes \cop)(\Phi)(\cop \otimes \id  \otimes \id)(\Phi)
\end{equation}
\begin{equation}
\label{eq: normal}
(\id \otimes \eps \otimes \id)(\Phi) = 1 \otimes 1.
\end{equation}

The next lemma applies the axioms above: 
\begin{lemma}
In any quasi-bialgebra $$\eps(X^1) X^2 \otimes X^3 = 1 \otimes 1 = X^1 \otimes X^2 \eps(X^3).$$
\end{lemma}
\begin{proof}
The known trick is to apply the counit to various equivalent forms of Eq.~(\ref{eq: 3-cocycle}).
The left equation follows from
$$ \eps (X^1) X^2 \otimes X^3  = Y^1 {Z^1}\1 x^1 y^1 \eps(Y^2 {Z^1}\2 x^2 {y^2}\1) {Y^3}\1 Z^2 x^3 {y^2}\2 
\otimes {Y^3}\2 Z^3 y^3  $$
$$ = Y^1 \eps(Y^2) Z^1 x^1 \eps(x^2) y^1 {Y^3}\1 Z^2 x^3 y^2 \otimes {Y^3}\2 Z^3 y^3 = 1 \o 1 $$
by two applications of Eq.~(\ref{eq: normal}) and finally $\Phi \Phi^{-1} = 1 \otimes 1 \otimes 1$. 
The right equation is similarly established.  
\end{proof}

These axioms mean  that the category of left or right modules over $H$ form a  non-strict tensor category
where multiplication by  Drinfel'd's associator $\Phi$ provides a natural isomorphism between triple tensor products of modules
(Eq.~(\ref{eq: quasi-coassoc})), all possible associations of the same tensor product of modules are isomorphic from
a commutative pentagon diagram 
(Eq.~(\ref{eq: 3-cocycle})) and the unit module $k$ with $H$-module structure induced from the augmentation $\eps$
is cancellable up to natural isomorphism in the middle position of a triple tensor product (Eq.~(\ref{eq: normal})).
A bialgebra is of course of quasi-bialgebra where $\Phi = 1 \otimes 1 \otimes 1$.  
Unlike bialgebra, the notion of quasi-bialgebra is stable under twisting of the coproduct $\cop \leadsto \mathcal{F} \cop
\mathcal{F}^{-1}$ \cite{Kas} where $\mathcal{F} \in H \otimes H$.  

A quasi-bialgebra $H$ is called a \textit{quasi-Hopf algebra} if there is an anti-automorphism $S: H \to H$ (called
an \textit{antipode}) with elements
$\alpha, \beta \in H$ such that for all $a \in H$:
\begin{eqnarray}
\label{eq: alpha}
S(a\1) \alpha a\2  & = & \eps(a) \alpha \\
\label{eq: beta}
a\1 \beta S(a\2)   & = & \eps(a) \beta \\
\label{eq: phi-beta-alpha}
X^1 \beta S(X^2) \alpha X^3 &= & 1 \\ 
 \label{eq: phi-inverse-alpha-beta}
S(x^1) \alpha x^2 \beta S(x^3) & = & 1
\end{eqnarray}

In $k$ we cannot assume a cancellation law, but it follows from Eqs.~(\ref{eq: normal}) and~(\ref{eq: phi-inverse-alpha-beta})
 that $\eps(\alpha) \eps(\beta) = 1$,
and then from Eq.~(\ref{eq: alpha}) or~(\ref{eq: beta}) that $\eps \circ S = \eps$. Since $\eps(\alpha)$ and
$\eps(\beta)$ are inverses of one another in $k$, we may rescale $\alpha$
and $\beta$ so that $\eps(\alpha) = 1 = \eps(\beta)$.  A Hopf algebra is of course a quasi-Hopf algebra where
$\alpha = \beta = 1$ and $\Phi = 1 \otimes 1 \otimes 1 $. 
 However, unlike for a Hopf algebra, the antipode
of a quasi-Hopf algebra is only unique up to inner automorphism of $H$:
 given another antipode $\overline{S}$, it is $S$ composed with an inner automorphism
with unit $u \in H$ where the transformation
is $\alpha \leadsto u\alpha$, $\beta \leadsto \beta u^{-1}$
and $S(a) \leadsto \overline{S} = u S(a)u^{-1}$ \cite{D}.  
The antipode $S$ also differs in general from a Hopf algebra antipode by
being only an anti-coalgebra automorphism up to a twist \cite{BC,D,HN}.

Four elements are introduced in order to generalize Hopf algebra formulae of the
type $a\1 \o a\2 S(a\3) = a \o 1$ to the quasi-Hopf setting.  They are the following
elements in $H \o H$:
\begin{eqnarray}
\label{eq: cue}
q_R := X^1 \o S^{-1}(\alpha X^3) X^2  & & q_L := S(x^1)\alpha x^2 \o x^3 \\
\label{eq: pee}
p_R := x^1 \o x^2 \beta S(x^3) & & p_L := X^2 S^{-1}(X^1 \beta) \o X^3  
\end{eqnarray}
Again briefly denote $q_R = q^1_R \o q^2_R$, etc.\ by suppressing the summation symbol and
indices.  
The formulae they facilitate are the following for each $a \in H$:
\begin{eqnarray}
\label{eq: cue-are}
q^1_R a_{(1,1)} \o S^{-1}(a\2) q^2_R a_{(1,2)}  & = & a q^1_R \o q^2_R \\ \label{eq: pee-are} 
a_{(1,1)} p_R^1 \o a_{(1,2)} p_R^2 S(a\2) & = & p_R^1 a \o p_R^2  \\
S(a\1) q_L^1 a_{(2,1)} \o q_L^2 a_{(2,2)} = q_L^1 \o a q_L^2 \\
\label{eq: pee-ell}
a_{(2,1)} p_L^1  S^{-1}(a\1) \o a_{(2,2)} p_L^2 = p_L^1 \o p_L^2 a 
\end{eqnarray}

For example, Eq.~(\ref{eq: cue-are}) follows from  variants of Eqs.~(\ref{eq: quasi-coassoc})
and~(\ref{eq: alpha}):
$$ LHS = X^1 a_{(1,1)} \o S^{-1}(a\2) S^{-1}(\alpha X^3)X^2 a_{(1,2)} $$
$$ = a\1 X^1 \o S^{-1}(\alpha
a_{(2,2)} X^3) a_{(2,1)} X^2 = RHS. $$

The element pairs $q_R$, $p_R$ and $q_L$, $p_L$ satisfy two equations each below:
\begin{eqnarray}
\label{eq: first}
\cop(q^1_R) p_R (1 \o S(q^2_R)) & = & 1 \o 1 \\
\label{eq: second}
(1 \o S^{-1}(p^2_R))q_R \cop(p^1_R) & = & 1 \o 1 \\
\cop(q^2_L) p_L (S^{-1}(q^1_L) \o 1) & = & 1 \o 1 \\
\label{eq: fourth}
(S(p^1_L) \o 1) q_L \cop(p^2_L) & = & 1 \o 1. 
\end{eqnarray}  
 
For example, Eq.~(\ref{eq: first}) follows from Eq.~(\ref{eq: normal}), the lemma, and 
Eq.~(\ref{eq: phi-beta-alpha}):
$$ LHS = {X^1}\1 x^1 \o {X^1}\2 x^2 \beta S(x^3) S(X^2) \alpha X^3 $$
$$ = x^1 X^1 \o 
x^2 Y^1 {X^2}\1 \beta S({x^3}\1 Y^2 {X^2}\2)\alpha {x^3}\2 Y^3 X^3  $$
$$ 
= 1 \o Y^1 \beta S(Y^2) \alpha Y^3 = 1 \o 1.$$

\subsection{A direct proof that $H$ is a Frobenius algebra}
In this subsection we show that the quasi-Hopf algebra $H$
is a Frobenius algebra by sketching the method in \cite{BC}.  
Suppose $\{ a_i \}_{i=1}^n$ and $\{ f^i \}^n_{i=1}$ are dual or projective bases in
$H$ and $H^*$, respectively: i.e., they satisfy only $\sum_{i=1}^n f^i(a) a_i = a$
for all $a \in H$, but not necessarily $f_i(a_j) = \delta_{ij}$.  However, in common with
dual bases for finite dimensional algebras over fields, we have: 
\begin{lemma}
In any f.g.\ projective $k$-algebra $H$, we have for every $x \in H$
\begin{equation}
\sum_i f^i \leftharpoonup x \o a_i = \sum_i f^i \o x a_i
\end{equation}
\begin{equation}
\sum_i x \rightharpoonup f^i \o a_i = \sum_i f^i \o a_i x.
\end{equation}
\end{lemma}
\begin{proof}
The proof is a brief calculation like in the case of a ground field once we see clearly
that $\sum_i \phi_i \o c_i = \sum_j \rho_j \o b_j $ in $H^* \o H$ if and only if
for each $a \in H, \eta \in H^*$ we have $\sum_i \bra \phi_i | a \ket \bra \eta | c_i \ket
= \sum_j \bra \rho_j | a \ket \bra \eta | b_j \ket $ by applying projective bases.  
\end{proof}

As we know from Hopf algebra theory, integrals are of interest in questions of
semisimplicity, or failing that, Frobenius/symmetric algebra properties.  

\begin{definition}
A left {\em integral} is an element $t \in H$ satisfying $at = \eps(a) t$ for all $a \in H$.  A
right integral $r \in H$ similarly satisfies $ra = \eps(a) r$ for $a \in H$. 
Denote the space of left integrals by $\int^{\ell}_H$ and right integrals by $\int^r_H$. 
\end{definition}
 
Following \cite{PO} and Van Daele one shows the existence of integrals in $H$ by defining
a projection $P : H \to \int^{\ell}_H$: ($h \in H)$

\begin{eqnarray}
\label{eq: Van Daele}
P(h) & = & \sum_{i=1}^n \bra f^i | \beta S^2(q^2_R {a_i}\2)h \ket \, q^1_R {a_i}\1 \\ 
& = & \sum_{i=1}^n \bra f^i | \beta S( S(X^2 {a_i}\2 ) \alpha X^3) h \ket \, X^1 {a_i}\1
\end{eqnarray}

We check that $P(h)$ is a left integral for each $h \in H$ using Eq.~(\ref{eq: cue-are}):  

$$ aP(h) =  \sum_i \bra f^i | \beta S^2(q^2 {a_i}\2) h \ket \, aq^1 {a_i}\1  $$
$$ = \sum_j \bra f^i | \beta S^2(S^{-1}(a\2)q^2 a_{(1,2)} {a_i}\2)h \ket \, q^1 a_{(1,1)} {a_i}\1 $$
$$ = \sum_{i,j} \bra f^i | a\1 a_j \ket \bra f^j | \beta S(a\2) S^2(q^2 {a_i}\2) h \ket \, q^1
 {a_i}\1 $$
 $$
 =  \sum_i \bra f^i | a\1 \beta S(a\2) S^2(q^2 {a_i}\2) h \ket\, q^1 {a_i}\1 = \eps(a) P(h)  $$
where $q_R = q^1 \o q^2$. 
In the third equation, we 
use $\cop(\sum_i a_i \bra f^i | xy \ket) =  x\1 y\1 \o x\2 y\2$.
The existence of a nonzero integral
 will follow if one of $b_j := P(a_j)$ ($j = 1,\ldots,n$) is nonzero:
using the antipode axioms again, we note 
$$ \sum_j \bra f^j | S(b_j \beta) \ket = \sum_{i,j} \bra 
f^i | \beta S(S(X^2 {a_i}\2) \alpha X^3 ) a_j \ket \bra f^j | S(X^1 {a_i}\1 \beta) \ket  $$
$$ = \bra f^i | \beta S(X^3) S(\alpha) S^2(X^2) S^2({a_i}\2) S(\beta) S({a_i}\1) S(X^1) \ket $$
$$  
= \bra \eps | \beta S(X^3) S(\alpha) S^2(X^2) S(\beta) S(X^1) \ket = 1, $$
from which the claim follows.

By using the elements $q_R, p_R$ defined in Eqs.~(\ref{eq: cue-are}) and~(\ref{eq: pee-are}) above, let 
\begin{equation}
\label{eq: undercop}
\underline{\cop}(h) := q^1_R h\1 p^1_R \o q^2_R h\2 p^2_R := h_{\underline{(1)}} \o h_{\underline{(2)}}
\end{equation}
and for $f \in H^*$ let $f \rightarrow h :=  f(h_{\underline{(2)}}) h_{\underline{(1)}}$.   
A computation exactly like that in \cite[2.2]{BC} establishes that over a commutative ground ring:  

\begin{theorem}[Bulacu-Caenepeel]
The mapping $\Theta: \int^{\ell}_H \o H^* \to H$ given by $$\Theta(t \o h^*) = h^*(S(t_{\underline{(2)}})) t_{\underline{(1)}}$$
is an isomorphism with respect to the natural left $H$-modules ${}_HH$ and ${}_HH^*$.  Its inverse is given by
$$ \Theta^{-1}(h) = \sum_{i=1}^n P(a_ih) \o f^i. $$
\end{theorem}

As a consequence, $\int^{\ell}_H$ is one-dimensional if $k$ is a ground field, since $\dim H = \dim H^*$.  Otherwise, all we can say is that 
$\int^{\ell}_H$ has constant rank 1 with respect to localizations at any prime ideal in $k$.  
If the Picard group of $k$ is trivial, e.g.\
when $k$ is a local, semilocal or polynomial ring, 
 this
will mean that $\int^{\ell}_H$ is free of rank one.  Somewhat more generally, 

\begin{definition}
We refer to a quasi-Hopf algebra $H$ with $\int^{\ell}_H$ free of rank one as a QFH-algebra. 
\end{definition}

We propose to also call a QFH-algebra, a ``quasi-Hopf-Frobenius algebra'' with the reverse FH symbolism 
as in \cite{Par72}.  

Recall that a Frobenius $k$-algebra $H$ is finitely generated projective over $k$ with
an isomorphism of the natural modules ${}_H H \cong {}_H H^*$ (equivalently $H_H \cong {H^*}_H$).
\cite{KS}. (A word of caution that such a Frobenius $k$-algebra is a quasi-Frobenius ring if $k$ is a quasi-Frobenius
commutative ring, but not for general $k$.)

\begin{cor}
A QFH-algebra $H$ is a Frobenius algebra.
\end{cor}
\begin{proof}
Let $t$ denote a nonzero left integral freely generating $\int^{\ell}_H$. 
Then $\theta: {}_HH^* \to {}_HH$ defined by $h^* \mapsto h^* \circ S \rightarrow t$ is
an isomorphism by theorem.     
\end{proof}

$\theta: H^* \to H$ is known as a \textit{Frobenius isomorphism}, from which a \textit{Frobenius coordinate system} for $H$ may be derived
via \textit{Frobenius homorphism} $\lambda := \theta^{-1}(1) \in H^*$, with dual bases $\{ b_i := \theta(f^i)\}$, $\{ a_i \}$.  From
the equation $\sum_i a_i f^i(a) = a$, all $a \in H$ we obtain the equation,
\begin{equation}
\label{eq: Frobenius1}
\sum_{i=1}^n \lambda(ab_i)a_i = a
\end{equation}
The symmetrical equation
 \begin{equation}
\label{eq: Frobenius2}
\sum_i b_i \lambda(a_i a) = a
\end{equation}
 for any $a \in H$ follows by noting $\theta^{-1}(a - \sum_i b_i \lambda(a_i a))$
on any $x \in H$ is zero. From either of these equations, we see that the Frobenius homomorphism $\lambda$ is
a \textit{nondegenerate functional} (i.e., $\lambda(Hx) = 0$ or $\lambda(xH) = 0$ implies $x = 0$).
We then see that $\lambda(ax) = \lambda(x \eta(a))$, or equivalently
$$\lambda \leftharpoonup a = \eta(a) \rightharpoonup \lambda,$$ defines an automorphism $\eta$ of $H$, 
which has the alternative definition 
\begin{equation}
\label{eq: Naka}
\eta(a) =  \sum_i b_i \lambda(a a_i) 
\end{equation} 
called the \textit{Nakayama automorphism}.

As an aside, we  point out that, like a symmetric algebra, a Frobenius algebra $H$ 
satisfies the bimodule isomorphism ${}_HH_H \cong {}_HH^*_{{\eta}^{-1}}$ 
where  we employ the notation for the 
obvious twist or pullback of module structure by an automorphism (which we will see later
in connection with $\beta$-Frobenius extensions). Recall that $H$ is 
a \textit{symmetric algebra} if $\eta$ is inner; equivalently, $H \cong H^*$ as $H$-bimodules.


\section{Separability and quasi-Hopf algebras}

We continue with the notation developed for a QFH-algebra in the last section and subsection.  

Recall that a separable $k$-algebra $A$ is characterized by having a separability element $e$
in $A \o A$ (or idempotent when
viewed  in $A^e := A \o A^{\rm op}$). Again suppressing summation and indices, we write $e = e^1 \o e^2$
and such a (nonunique) $e$ is characterized by $e^1 e^2 = 1$ and a \textit{Casimir condition} $ae^1 \o e^2 = e^1 \o e^2a$ for all $a \in H$.
With this element we may also characterize $A$ by all $k$-split exact sequences  of $A$-modules are in fact $A$-split exact
(using the Maschke technique of applying the separability element
to the argument and value of a function).  Over a commutative
ground ring $k$, a separable algebra $A$ is not necessarily semisimple;
however, if $k$ is semisimple, then $A$ is semisimple. 

 If $A$ is $k$-separable and f.g.\ projective, faithful over $k$,
the  Endo-Watanabe Theorem shows by complicated arguments (or ``big machinery'') that $A$ is a Frobenius algebra. 
 However, if $A$
is already known to be a symmetric algebra, we may apply a simple test to a Frobenius system $\phi: H \to k$
and $\{ x_i \}$, $\{ y_i \}$:
$A$ is $k$-separable if and ony if there is $a \in A$ such that $\sum_i x_iay_i = 1$.  This is proven by 
using ideas from the proof of Lemma~\ref{lem-preRadford} below as well as noting that 
the equations~(\ref{eq: Frobenius1}) and~(\ref{eq: Frobenius2}) imply
that $\sum_i x_i \o y_i$  satisfies a Casimir condition. 

The next theorem is Panaite's when $k$ is a field.  

\begin{theorem} 
\label{thm-sep}
A quasi-Hopf algebra $H$ is separable over its commutative ground ring $k$ if and only if there is a normalized
left or right integral in $H$. 
\end{theorem}

\begin{proof}
($\Rightarrow$)  Let $K = \ker \eps$, a two-sided ideal in $H$.  The counit $\eps$ induces a $k$-split exact sequence
$$ 0 \longrightarrow {}_HK \longrightarrow {}_HH \stackrel{\eps}{\longrightarrow} {}_H k \longrightarrow 0$$
where the last nonzero module is induced by the augmentation $\eps$.  By $k$-separability of $H$, this short exact
sequence is split over $H$, whence there is a left $H$-homomorphism $t' : k \to H$ such that $\eps t' = \id_k$.  
But $t := t'(1)$ is then a normalized left integral, since $ht = t'(h \cdot 1) = \eps(h)t$, all $h \in H$
 and $\eps(t) = 1$. (I.e., ${\rm Hom}_{H-}\, (k,H) \cong \int^{\ell}_H$ as Wisbauer has observed.) This argument
may be repeated with right $H$-modules to establish a normalized right integral (also without the presence of antipode).

($\Leftarrow$) Given a normalized left integral $t$ or right integral $r$, any of the following four are separability elements:
\begin{eqnarray}
e_{1,2} & = & S(r\1 p^1) \o \alpha r\2 p^2 \ \ \ \ \ \ \ (p = p_L \ {\rm or} \ p_R) \\
e_{3,4} & = & q^1 t\1 \beta \o S(q^2 t\2)  \ \ \ \ \ \ \ (q = q_L \ {\rm or} \ q_R )
\end{eqnarray}
For example, $e_1$ is a separability element, since
 $S(p^1_L)S(r\1)\alpha r\2 p^2_L = S(p^1_L) \alpha p^2_L =1$ by applying $\id \o \eps$ to
Eq.~(\ref{eq: fourth}) and noting $q^1_L \eps(q_L^2) = \alpha$ from Eq.~(\ref{eq: cue}). The Casimir condition
follows from combining Eq.~(\ref{eq: pee-ell}) with the trivial observation $\cop(ra) = 
\cop(r) \eps(a)$, all $a \in H$. ($e_4$ satisfies the  Casimir condition  and $\beta e_4 = 
t\1 \beta \o S(t\2)$ by \cite[2.1]{BC}.)
\end{proof}

In fact, if H is a separable algebra, left and right normalized integrals coincide in what we call the \textit{Haar
 integral}, i.e., the algebra $H$ is \textit{unimodular}, for the following  
reason which holds for an augmented Frobenius algebra and 
 makes use of the \textit{modular augmentation}
$\mu : H \to k$ defined by
\begin{equation}
\label{eq: modular}
ta = \mu(a)t.
\end{equation}
\begin{prop}
The counit, Nakayama automorphism and modular augmentation in an augmented Frobenius algebra satisfy
\begin{equation}
\eps = \mu \circ \eta.
\end{equation}
Consequently, an augmented symmetric algebra is unimodular.
\end{prop}
\begin{proof}
Given counit $\eps$ in an augmented Frobenius algebra $A$, there is a left integral $t \in A$
and Frobenius homomorpism $\phi \in A^*$ such that $\phi(t) = 1$ \cite{KS}. 
Now define Nakayama automorphism $\eta$
and modular augmentation $\mu$ relative to $\phi$ and $t$ as above.  
 Whence $$ \mu =  \phi \leftharpoonup t= 
\eta(t) \rightharpoonup \phi   = (t \rightharpoonup \phi) \circ \eta^{-1} = \eps \circ \eta^{-1}, $$
since $\phi \circ \eta = \phi$.  

If $A$ is a symmetric algebra, $\eta(a) = u a u^{-1}$ for some unit $u \in A$ and $a \in A$.
 It follows that $\mu = \eps$, hence $A$ is
unimodular.
\end{proof}
 
Since a separable f.g.\ projective faithful $k$-algebra is a symmetric algebra, it follows
that a separable quasi-Hopf algebra $H$ is unimodular.    

In characteristic $p$, there is the phenomenon of strong separability in the sense of Kanzaki
\cite{KS2} which stands out as a strong form of separability.  In case $k$ is an algebraically
closed field, these are separable algebras all of whose simple modules have dimension relatively prime
to the characteristic of $k$. For general ground ring $k$, a separable $k$-algebra is
\textit{strongly separable} if $A = C \oplus [A,A]$ where $C$ is its center and $[A,A]$ is the $k$-span
of commutators $[a,b] = ab - ba$. Equivalently, $A$ is strongly separable if it has a symmetric
separability element, or even weaker,  an element $e \in A \o A$ such that $e^1 e^2 = 1$
and $e^1 a \o e^2 = e^1 \o a e^2$ for each $a \in A$ \cite{KS2}. 

 In \cite[4.1]{KS2}, the following criterion is given: a Frobenius algebra $A$ with 
Frobenius homomorphism $\phi: A \to k$ and dual bases $\{ x_i \}$, $\{ y_i \}$ (such
that $\sum_i x_i \phi(y_i a) = a$) is strongly separable if and only if 
$\sum_i y_i x_i$ is an invertible element $u$ in $A$.  Moreover, the Nakayama automorphism
is given by $\eta(a) = u^{-1} a u$ (and $A$ is naturally a symmetric algebra).  We
apply this criterion next to extend an old result of Larson for Hopf algebras \cite{KS2}.

\begin{theorem}
Suppose $H$ is a  $k$-separable quasi-Hopf algebra 
where $\beta   S(\alpha)  = 1$ and $S^2 = \id$ 
for some antipode $S$ on $H$. 
Then $H$ is strongly separable and $\lambda$ is a trace. 
\end{theorem}
\begin{proof}
Let $t$ be a Haar integral in $H$.  
Recalling dual bases $ \{ f^i \circ S \rightarrow t \}$, $\{ a_i \}$, and 
Frobenius homomorphism $\lambda$
from Section~2, we  show that 
$u = \sum_i a_i t_{\underline{(1)}}f^i(S(t_{\underline{(2)}}))=1$
below using the two lemmas directly following the proof: set
$p_R = p^1 \o p^2, q_R = q^1 \o q^2$, and note that 
\begin{eqnarray*}
S(\alpha) a_i \beta S(\alpha) (f^i \circ S \rightarrow t) & = & S(\alpha)S(p^2)S(t\2)S(q^2)\beta
S(\alpha)q^1 t\1 p^1 \\
&= & S(t\2 p^2 \alpha)S(\alpha)t\1 p^1 \\
& = & S(t\2)S(\alpha)t\1 = S(\alpha), 
\end{eqnarray*}
the last equation using \ref{eq: alpha}. 
Since $H$ is a f.g.\ projective $k$-algebra, $S(\alpha)$ is invertible, and the computation implies
that $u = 1$. Then $H$ is strongly separable with $\lambda$  a trace
since the Nakayama automorphism $\eta = \id$. 
\end{proof}

Note too the dual bases tensor is symmetric under the hypotheses.

\begin{lemma}
If $t \in \int^{\ell}_H$, then $q^1 t\1 \o S^{-1}(\beta) q^2 t\2 = t\1 \o t\2$.
\end{lemma}
\begin{proof}
Multiply Eq.~(\ref{eq: second}) from the right by $\cop(t)$, use the left integral property
and $\eps(p^1)p^2 = \eps(x^1)x^2 \beta S(x^3) = \beta$. 
\end{proof}

\begin{lemma}
If $r \in \int^r_H$, then $r\1 p^1 \o r\2 p^2 \alpha = r\1 \o r\2$.  
\end{lemma}
\begin{proof}
Multiply Eq.~(\ref{eq: first}) from the left by $\cop(r)$ and note $\eps(q^1)q^2 = S^{-1}(\alpha)$.
\end{proof}

From Theorem~\ref{thm-sep} it follows that we also have
\begin{eqnarray}
\beta q^1 t\1 \o S(q^2t\2) &= & t\1 \o S(t\2) \\
r\1 p^1 S^{-1}(\alpha) \o r\2 p^2 & = & r\1 \o r\2 .
\end{eqnarray}


\section{Frobenius coordinate systems and Radford's formula}

We begin this section with a basic lemma for a Frobenius algebra $A$ with anti-automorphism $S$
which establishes an archetypical result for Radford formulas for the fourth power of an antipode
in some quantum algebra.  It is stated in terms of a Nakayama automorphism which for a unimodular Hopf algebra is 
the square of $S$ or its inverse.  
We continue the conventions
begun above.  Let ${\rm Ad}_u$ denote conjugation by a unit $u$ where ${\rm Ad}_u(x) = uxu^{-1}$.  

\begin{lemma}[The Pre-Radford Formula]
\label{lem-preRadford}
If $A$ is a Frobenius algebra with Nakayama automorpism $\eta$ and anti-automorphism $S: A \to A$,
then there is invertible $d \in A$ such that
\begin{equation}
\label{eq: preRadford}
S\circ \eta \circ S^{-1} \circ \eta = {\rm Ad}_{d^{-1}}.
\end{equation}
\end{lemma}
\begin{proof}
Suppose $(\phi \in A^*, x_i, y_i, \eta)$ is the Frobenius system with Nakayama automorphism $\eta$ on $A$.
Let's recall (from any of several elaborative sources, e.g.\ \cite{KS,KS2,KS3}) 
that any other Frobenius system $(\psi \in A^*, u_j, v_j, \rho)$ only differs
from the first by an invertible element $d \in A$ (called the (right) \textit{derivative} $\frac{d\phi}{d\psi}$) such that 
\begin{eqnarray}
\psi & = & \phi \leftharpoonup d \\
\label{eq: dee}
d & = & \sum_i \psi(x_i)y_i \\
\sum_j u_j \o dv_j & = & \sum_i x_i \o  y_i \\
\eta^{-1} \circ \rho & = & d (-) d^{-1} 
\end{eqnarray}
The first equation follows from the fact that the $k$-dual $A^*$ is freely generated by each 
Frobenius homomorphism. The second follows from the first and Eq.~(\ref{eq: Frobenius1}).
The third equation follows the second and the Casimir condition for dual base tensors.
The fourth equation follows from the computation: ($a \in A$)
\begin{eqnarray}
\rho(a) &= & \sum_j u_j \psi(a v_j) \nonumber  \\
        & = & \sum_i x_i \phi(dad^{-1}y_i) \nonumber  \\
\label{eq: Naka transform}
        & = & \sum_i x_i \phi(y_i \eta(dad^{-1})) = \eta(dad^{-1}).
\end{eqnarray}

Next we claim that $S$ transforms a Frobenius system into another as follows: 

$$ (\phi, x_i,y_i, \eta) \leadsto \ (\phi \circ S^{-1},\, S(y_i),\, S(x_i),\, S \circ \eta^{-1} \circ S^{-1})$$

This is due to $$S^{-1}(a) = \sum_i \phi(S^{-1}(a) x_i)y_i = \sum_i y_i \phi \circ S^{-1}(S(x_i)a),$$
to which we apply $S$. We compute the Nakayama automorphism $\rho$ 
associated to $\psi  := \phi \circ S^{-1}$: ($a \in A$) 
\begin{eqnarray*}
\rho(a) & = & \sum_i S(y_i) \phi (S^{-1}(aS(x_i))) \\
        & = & S(\sum_i y_i \phi(x_i S^{-1}(a))) \\
        & = & S(\sum_i \phi (\eta^{-1} S^{-1}(a) x_i) y_i) = S(\eta^{-1}(S^{-1}(a))).
\end{eqnarray*}
Combining the existence of invertible $d \in A$ such that $\rho = \eta (d(-)d^{-1})$
with this last result in \cite{KS2}, we conclude that $S \circ \eta \circ S^{-1} \circ \eta = {\rm Ad}_{d^{-1}}$.   
\end{proof}

This lemma may be viewed as a key to understanding several Radford formulas for the fourth power
of antipodes on Hopf algebras, weak Hopf algebra, quasi-Hopf algebras and future
quantum Frobenius algebras. The Nakayama automorphism is often expressible in terms of the second
power of the antipode, whence the left-hand side of Eq.~(\ref{eq: preRadford}) will
involve the fourth power of antipode.  
 In support of 
this claim let us briefly consider the first two cases before we take up the third case 
in more detail later in this section.

 Let $H$ be a Hopf algebra, finite projective over a commutative ring $k$, 
 with right integral
$t \in H$ and right integral $f$ on $H$ such that $f(t) = 1$ (whence $f \leftharpoonup t = \eps$).
The conceptually brief proof of Radford's formula below is based on \cite{KS,KS3,K}.  
Since $S(a) = f(t\1 a) t\2$ satisfies $S(a\1)a\2 = \eps(a) 1_H$ for all $a \in H$, it defines the antipode $S: H \to H$
and it follows directly that a Frobenius system is given by 
$(f, \, S^{-1}(t\2), \, t\1)$.   The Nakayama automorphism of $f$ is given by
$$\alpha(x) = S^{-2}(x) \leftharpoonup m = S^{-2}(x \leftharpoonup m)$$
where $m: H \to k$ is the modular augmentation such that $at = m(a)t$ for all $a \in H$,
since $a \leftharpoonup f = f(a) 1_H$ and we apply $S^2$ to $\alpha(a) = S^{-1}(t\2)f(at\1)$. 
Consider now the anti-automorphism $S^{-1}$ on $H$.
With $f \circ S =  f \leftharpoonup d$ and $b$  the distinguished group-like element satisfying
 $\gamma f = \gamma(b)f$ for every $\gamma \in H^*$, we compute that both $f \circ S$ and $f \leftharpoonup b$
belong to the free rank one $k$-module of left integrals (since $b$ is grouplike with inverse $S(b)$) and assume the same 
value on the left integral $S^{-1}(t)$ since $f(S^{-1}(t)) = 1$ and $\eps(b) = 1$. 
Then $f \leftharpoonup b = f \leftharpoonup d$, whence $b = d$.  From the Lemma,
$S^{-1} \circ \alpha \circ S \circ \alpha = {\rm Ad}_{b^{-1}}$. Applying this to an $x \in H$ yields
$$ S^{-1}( S^{-2}(S^{-1}(x \leftharpoonup m)) \leftharpoonup m) = b^{-1}xb.$$ Since 
$S^{-1}(x \leftharpoonup m) = m^{-1} \rightharpoonup S^{-1}(x)$ for $x \in H$, this last equation
simplifies to $$m^{-1} \rightharpoonup S^{-4}(x) \leftharpoonup m = b^{-1}xb.$$  
A simplification yields Radford's formula  for the fourth power of the antipode:

\begin{equation}
\label{eq: Radford formula}
 S^4(x) = b(m^{-1} \rightharpoonup x \leftharpoonup m) b^{-1} 
\end{equation}

Consider next a special case of weak Hopf algebra $A$ in \cite[p. 423]{BNS}:
assume the existence of two-sided integrals that are nondegenerate, $h \in A$
and $h^* \in A^*$.  Define $a_L = h^* \rightharpoonup h$ and $a_R = S(a_L)$.
By \cite[Eq.\ (3.60)]{BNS}, $h^*$ is nondegenerate with Nakayama automorphism
$\theta_{h^*} = {\rm Ad}_{a_R} \circ S^2$.  By \cite[Lemma 3.20]{BNS},
$S^*(h^*) = h^*$.  Then the lemma above with derivative $d = 1$ shows
$$\theta_{h^*} = S \circ \theta^{-1}_{h^*} \circ S^{-1}$$
 from which it follows that ${\rm Ad}_{a_R} \circ S^2 = {\rm Ad}_{a_L} \circ S^{-2}$; whence
another way to see \cite[eq. (3.51)]{BNS}

\begin{equation}
S^4 = {\rm Ad}_{a^{-1}_R a_L}
\end{equation}

We return to our approach to a QFH-algebra $H$ via Frobenius system
$(\lambda, b_i, a_i, \eta)$ defined in Section~2.  We need a formula for
the Nakayama automorphism $\eta: H \to H$ in terms of $S$ such as in Hausser-Nill \cite[5.1]{HN}
for another Nakayama automorphism.
We will connect our approach via Bulacu-Caenepeel \cite{BC}
with the Hausser-Nill theory of left cointegrals and quasi-Hopf bimodules as follows.

We first briefly recall the notation 
$$\overline{\cop}(x) = x_{\overline{(1)}} \o x_{\overline{(2)}} := V\cop(x)U$$
for two invertible elements $U,V \in H \o H$ defined in \cite[(3.7)]{BC}
and \cite[3.12-3.13]{HN}.

\begin{lemma}
\label{lem-HNsystem}
   The Frobenius homomorphism $\psi := \lambda \circ S$ is a left cointegral
with dual bases tensor 
\begin{equation}
\label{eq: equality}
t_{\underline{(2)}} \o S^{-1}(t_{\underline{(1)}}) = 
S(r_{\overline{(1)}}) \o r_{\overline{(2)}}
\end{equation}
where $r = S^{-1}(t)$ is the right integral satisfying $\psi(r) = 1$. Its Nakayama automorphism
is given by 
\begin{equation}
\label{eq: HN-naka}
\rho(a) = S( S(a) \leftharpoonup \mu)
\end{equation}
for the modular augmentation $\mu: H \to k$.  
\end{lemma}

\begin{proof}
Given an augmented Frobenius algebra $(A, \eps)$ with Frobenius homomorphism $\psi$
(a free generator of $A^*$),
the left integrals are free of rank one; with left integral $t$ satisfying $\psi(t) = \eps$
(i.e., $t$ is a \textit{left norm}),
it is easy to show that $a \mapsto \psi(a)t$ is a projection onto the left integrals \cite{KS}. 
Moreover, for any $a \in A$ we note that $ta = t \mu(a)$ for another augmentation $\mu$
called the \textit{modular augmentation},
since $t$ freely generates $ta \in \int^{\ell}_A$ \cite{KS}.    

Next recall that the Frobenius homomorphism $\lambda$ was defined above via $\Phi^{-1}(1) =
\sum_i P(a_i) \o f^i$ where $P : H \to \int^{\ell}_H$ is another projection of $H$ onto
the left integrals and $\{ a_i \}$, $\{ f^i \}$ projective bases for $H$, $H^*$.  With $t\in \int^{\ell}_H$ such that $\lambda(t) = 1$, we see that
$\lambda = \sum_i \lambda(P(a_i)) f^i$. But for any $a \in H$, $P(a) = \sum_i P(a_i) f^i(a)$
from which it follows that
\begin{equation}
\label{eq: P-invariance}
\lambda(a) = \lambda(P(a))
\end{equation}
for each $a \in H$.  
(It follows that $P$ and the projection $a \mapsto \lambda(a)t$ are one and the same.)

We recall from Hausser-Nill that the space of left cointegrals in $H^*$ is
$\mathcal{L} = E(H^*)$ for a projection $E$ defined on $H^*$ in \cite[(3.3),(4.5)]{HN}.
The precise relationship between $E$ and $P$ is noted in \cite{BC} as  
follows:
\begin{equation}
\label{eq: BC}
\bra E(h^*) | h \ket = \bra h^* | S^{-1} P S (h) \ket
\end{equation}
To show $\psi = \lambda \circ S$ a left cointegral, it suffices to show $E(\psi) = \psi$
from the two equations~(\ref{eq: BC}) and~(\ref{eq: P-invariance}):
\begin{eqnarray*}
\bra E(\psi) | h \ket & = & \bra \lambda \circ S | S^{-1} P S(h) \ket \\
                      & = & \bra \lambda | PS(h) \ket \\ 
                      & = & \bra \lambda | S(h) \ket = \bra \psi | h \ket
\end{eqnarray*}
for each $h \in H$. 

The results in \cite{HN} are valid for QFH-algebra $H$ because they require
only the Drinfeld calculus introduced in the preliminaries, as well as
 for the following reason.  
Since $\int^{\ell}_H$ and $\int^r_H$ are isomorphic as $k$-modules under $S$,
both are free of rank one.   But $\int^r_H$ and $\mathcal{L}$ are nondegenerately
paired by \cite[Lemma 4.4]{HN}, which shows that $\mathcal{L}$ is also free
of rank one.  

Next we note that $\psi(r) = \lambda S(S^{-1}(t)) = \lambda(t) = 1$ by choice of $t$, and 
that \cite[Prop.\ 5.5]{HN} shows that $S(r_{\overline{(1)}}) \o r_{\overline{(2)}}$
as the dual bases tensor for the unique left cointegral $\psi$ such that $\psi(r) = 1$. 

From the lemma, $\psi$ has dual bases tensor $\sum_i S^{-1}(a_i) \otimes S^{-1}(b_i)$.  
Recalling that
$b_i = t_{\underline{(1)}} f^i(S(t_{\underline{(2)}}))$, it follows that
the dual bases tensor is the left-hand side of Eq.~(\ref{eq: equality}).

Finally Hausser and Nill show the Nakayama automorphism for $\psi$ in \cite[Lemma 5.1]{HN}
to be $\rho = S \circ S_{\mu}$ where $S_{\mu}(a) := S(a) \leftharpoonup \mu$
and $\mu$ is the modular
augmentation.
\end{proof}

From this lemma and Eq.~(\ref{eq: Naka transform}) the formula 
for the inverse of the Nakayama automorphism $\eta$  introduced in Eq.~(\ref{eq: Naka}) 
is seen to be 
\begin{equation}
\eta^{-1}(a) = S^2 (a \leftharpoonup \mu).
\end{equation}
It follows that the outer automorphism coset of the Nakayama automorphism does not change
upon changing antipode for $H$ (cf.\ preliminaries).

Let $u := \psi(r_{\overline{(1)}}) r_{\overline{(2)}}$, the comodulus or distinguished group-like
element in $H$. The lemma just proven has two consequences.
 
\begin{theorem}
\label{th-HN}
A QFH-algebra $H$ has antipode $S$ satisfying the Hausser-Nill equation:
\begin{equation}
\label{eq: HN}
S^2 \circ S^2_{\mu} = {\rm Ad}_{u^{-1}}
\end{equation}
\end{theorem}
\begin{proof}
We apply Lemma~\ref{lem-preRadford} to the Frobenius homomorphism
 $\psi = \lambda \circ S$ with Nakayama automorphism $\rho = S \circ S_{\mu}$.
We first compute $d = u$ from Eq.~(\ref{eq: dee}) in transforming from $\psi$
into $\psi \circ S^{-1} = \lambda$. Then
$$ S \rho S^{-1} \rho = S^2 S_{\mu}^2 = {\rm Ad}_{d^{-1}}. \qed$$
\renewcommand{\qed}{}\end{proof}

For an application of the ideas in Lemma~\ref{lem-HNsystem}
to the unimodularity problem for the Drinfel'd double $D(H)$, see Bulacu and Torrecillas \cite{BT}. 

In this paper, a \textit{quasi-Hopf subalgebra} $K \subseteq H$ is a 
$k$-subalgebra such that $K$ is a pure
$k$-submodule of $H$ \cite{Lam} for which $\cop(K) \subseteq K \o K$
and $K$ has its own associator $\Phi_K \in K \o K \o K$ \cite{Sch}; in addition, we assume
of our quasi-Hopf subalgebra 
that $H$ has an antipode $S$ stabilizing $K$ ($S(K) \subseteq K$) and that there are 
elements $\alpha_K, \beta_K$ in $K$ which together with  $S$ 
satisfy the axioms~(\ref{eq: alpha})-(\ref{eq: phi-inverse-alpha-beta}). 
It follows from some pure module theory
that $K$ is f.g.\ projective as a $k$-module \cite{Lam}. 
  A QFH-subalgebra pair $K \subseteq H$
is a quasi-Hopf subalgebra where both are QFH-algebras (so both are Frobenius algebras). 

Recall that a subring pair $R \supseteq S$ is called a $\beta$-Frobenius extension
if $\beta: S \to S$ is a ring automorphism, the natural module $R_S$ is f.g.\ projective 
with $\Hom (R_S, S_S) \cong R$
as $S$-$R$-bimodules where all module actions are the natural ones except the left $S$-module
structure on $S$ is the pullback module ${}_{\beta}S$ under the mapping $\beta: S \to S$. 
There are close connections explored by Kasch, Nakayama-Tzuzuku and Hirata
 between the module categories of $R$ and $S$ in such an extension.   
The next theorem generalizes facts obtained by Oberst-Schneider
and Fischman-Montgomery-Schneider \cite{FMS}.

\begin{theorem}
\label{th-beta}
A QFH-subalgebra pair $K \subseteq H$ forms a $\beta$-Frobenius extension where
$\beta$ is the relative Nakayama automorphism $\rho_K^{-1} \circ \rho_H$.  
\end{theorem}

\begin{proof}
The proof will follow from  \cite[Prop.\ 5.1]{KS} where we 
recall that a Frobenius subalgebra pair such as $K \subseteq H$ is $\beta$-Frobenius
with $\beta$ the relative Nakayama automorphism of $K$ if
two conditions are met:
\begin{enumerate}
\item $H_K$ is f.g.\ projective;
\item the Nakayama automorphism $\rho_H$ of $H$ stabilizes $K$: $\rho_H(K) \subseteq K$. 
\end{enumerate}
The second condition is met by $\rho_H(a) = S (S(a) \leftharpoonup \mu)$ ($a \in H$) since
$K$ is stable under $S$ and $\cop$.  The first condition follows from Schauenburg's freeness
theorem for quasi-Hopf subalgebras over ground fields \cite[3.2]{Sch} and two lemmas below adapting
this to commutative ground rings via localization at maximal ideals.
\end{proof} 

A formula for a Frobenius homomorphism $F: H \to K$ from \cite[5.1]{KS} is given by
\begin{equation}
\label{eq: Frob. homo}
F(a) = \psi (a \Lambda_{\underline{(2)}}) S^{-1}(\Lambda_{\underline{(1)}})
\end{equation}
where $\psi$ is the Frobenius homomorphism for $H$ in Lemma~\ref{lem-HNsystem}
and free generator $\Lambda \in \int^{\ell}_K$.  

Recall that projective modules over local rings are free, so that the next lemma is valid for
QFH-algebra over a local ring. 

\begin{lemma}
If $H$ is a finitely generated  free quasi-Hopf algebra over a local ring $k$
with $K$ a quasi-Hopf subalgebra, then the natural modules $H_K$ and ${}_KH$ are free.
\label{lemma-topfjell}
\end{lemma}
\begin{proof}
It will suffice by symmetry to prove that $H_K$ is free. 
First note that $H_K$ is finitely generated since $H_k$ is. 
If $ \mathcal{M}$ is  the maximal ideal of $k$, then the finite dimensional
quasi-Hopf algebra $\overline{H} := H/\mathcal{M}H$ is  free over the quasi-Hopf subalgebra 
$\overline{K} := K/\mathcal{M}K$
by purity and the freeness theorem in \cite{Sch}. Suppose $\theta: \overline{K}^n \stackrel{\cong}{\rightarrow}
\overline{H}$ is a $\overline{K}$-linear isomorphism.  Since $K$
is finitely generated over $k$, 
$\mathcal{M}K$ is contained in the radical of $K$. 
Now  $\theta$ lifts
to a right $K$-homomorphism $K^n \rightarrow H$ with respect to the natural projections
$H \rightarrow \overline{H}$ and $K^n \rightarrow \overline{K}^n$. 
By Nakayama's lemma, the homomorphism $K^n \rightarrow H$ is  epi (cf.\ \cite{Sil}).  
Since $H_k$ is finite projective, $\tau$ is a $k$-split epi, 
which is bijective
by Nakayama's lemma  applied to the underlying $k$-modules. 
Hence, $H_K$ is free of finite rank. 
\end{proof}

Over a non-connected ring $k = k_1 \times k_2$,
it is easy to construct examples of   (quasi-)Hopf subalgebra pairs
\[
K := k[H_1 \times H_2] \subseteq H := k[G_1 \times G_2]
\]
 where $G_1 > H_1$,
$G_2 > H_2$ are  subgroup pairs of finite groups and $H_K$ is not free 
(by counting dimensions on either side of $H \cong K^n$). 
  The next lemma follows from the previous one. 

\begin{lemma}
If $H$ is a $k$-finite projective, quasi-Hopf algebra and $K$ is a 
quasi-Hopf subalgebra of $H$, then the natural modules $H_K$ and ${}_KH$
are finite projective. 
\label{lemma-fjelltop}
\end{lemma}
\begin{proof}
 First note that $H_K$ is finitely generated.
If $k$ is a commutative ground ring, $Q \rightarrow P$ is an epimorphism
of $K$-modules, then it will suffice to show
that the induced map $\Psi: {\rm Hom}_{-K}(H, Q) \rightarrow {\rm Hom}_{-K}(H,P)$ is
epi too. Localizing at a maximal ideal
$\mathcal{M}$ in $k$, we obtain local ring $k_{\mathcal{M}}$, modules
over $k$ localized over $\mathcal{M}$, and a homomorphism denoted by $\Psi_{\mathcal M}$ as follows. 
By adapting a standard argument such as in \cite{Sil}, we note that for every module $M_K$
\begin{equation}
\Hom_K(H_K,M_K)_{\mathcal{M}} \cong \Hom_{-K_{\mathcal{M}}}(H_{\mathcal{M}},
M_{\mathcal{M}})
\label{eq:isom}
\end{equation}
since $H_k$ is finite projective.  Then $\Psi_{\mathcal M}$
maps 
$$ \Hom_{-K_{\mathcal{M}}}(H_{\mathcal{M}},
Q_{\mathcal{M}}) \rightarrow \Hom_{-K_{\mathcal{M}}}(H_{\mathcal{M}},
P_{\mathcal{M}}).$$

By Lemma~\ref{lemma-topfjell}, $H_{\mathcal{M}}$ is free over
$K_{\mathcal{M}}$.  It follows that $\Psi_{\mathcal{M}}$ is epi
for each maximal ideal $\mathcal{M}$, whence $\Psi$ is epi. 
\end{proof}



\begin{thebibliography}{XXXXXX}
\begin{small}
\bibitem{BNS}{G.~B\"ohm, F.~Nill and K.~Szlach\'anyi,
Weak Hopf algebras, I. Integral theory and $C^*$-structure,
\textit{J. Algebra} \textbf{221} (1999), 385-438.}
\bibitem{B}{D.~Bulacu,
Relative Hopf modules for (dual) quasi-Hopf algebras,
\textit{J.\ Algebra} \textbf{229} (2000), 632--659.}
\bibitem{BC}{D. Bulacu and S. Caenepeel,
Integrals for (dual) quasi-Hopf algebras. Applications.  \textit{J.\ Algebra},
\textbf{266} (2003), 552--583.} 
\bibitem{BT}{D.~Bulacu and B.~Torrecillas,
Factorizable quasi-Hopf algebras-Applications. \textit{J.\ Pure Appl.\ Alg.},
to appear.}  \texttt{math.QA/0312076}.
\bibitem{D}{V.G. Drinfel'd,
Quasi-Hopf algebras, \textit{Leningrad Math.\ J.} \textbf{1} (1990), 1419--1457.}
\bibitem{EW}{S. Endo and Y. Watanabe, 
On separable algebras over a commutative ring,
{\it Osaka J. Math.} {\bf 4} (1967), 233--242.}
\bibitem{FMS}{D. Fischman, S. Montgomery, and H.-J. Schneider,
Frobenius extensions of subalgebras of Hopf algebras,
{\it Trans.\ Amer.\ Math.\ Soc.\ } {\bf 349} (1997), 4857--4895.}
\bibitem{HN}{F. Hausser and F. Nill,
Integral theory for quasi-Hopf algebras, preprint, \texttt{math.QA/9904164}.}
\bibitem{K}{L.~Kadison,
\textit{New Examples of Frobenius Extensions}, University Lecture Series \textbf{14},
A.M.S., Providence, 1999. Update, 6 pp: www.ams.org/bookpages.}
\bibitem{KS2}{L. Kadison and A.A. Stolin,
Separability and Hopf algebras, in: \textit{Algebra and its Applications} 
(Athens, Ohio, March 1999), eds.\ Huynh, Jain
and Lopez-Permouth, Contemporary Math.\ \textbf{259}, A.M.S.,
Providence, (2000), 279--298.}
\bibitem{KS}{L. Kadison and A.A. Stolin,
An approach to Hopf algebras via Frobenius coordinates, \textit{Beitr\"{a}ge
z. Alg.\ u. Geometrie} \textbf{42} (2001), 359--384.}
\bibitem{KS3}{L. Kadison and A.A. Stolin,
An approach to Hopf algebras via Frobenius coordinates II, \textit{J.\ Pure Appl.\ Alg.}
\textbf{176} (2002), 127--152.}
\bibitem{Kas}{C. Kassel,
\textit{Quantum Groups}, Grad.\ Texts Math.\ {\bf 155}, Springer, New York, 
1995.} 
\bibitem{Lam}{T.Y. Lam,
{\it Lectures on Modules and Rings}, Grad.\ Texts Math.\ {\bf 189}, Springer-Verlag,
Heidelberg-Berlin-New York, 1999.}
\bibitem{PO}{F. Panaite and F. Van Oystaeyen,
Existence of integrals for finite dimensional quasi-Hopf algebras, \textit{Bull.\ Belg.\
Math.\ Soc.-Simon Stevin} \textbf{7} (2000), 261--264.}
\bibitem{Par72}{B. Pareigis,
On the cohomology of modules over Hopf algebras, {\it J. Algebra} {\bf 22}
(1972), 161-182.}
\bibitem{Sch}{P. Schauenburg,
A  quasi-Hopf algebra freeness theorem, preprint, \textit{Proc.\ Amer.\ Math.\ Soc.} \textbf{132}
(2004), 965--972.}
\bibitem{Sil}{J.R. Silvester,
\textit{Introduction to Algebraic K-theory}, Chapman and Hall, London, 1981.}
\end{small}
\end{thebibliography}
\end{document}